\documentstyle[12pt,amscd]{amsart}
\newtheorem{Theorem}{Theorem}[section]
\newtheorem{Lemma}[Theorem]{Lemma}
\newtheorem{Corollary}[Theorem]{Corollary}
\newtheorem{Proposition}[Theorem]{Proposition}
\newtheorem{Example}[Theorem]{Example}
\newtheorem{Remark}[Theorem]{Remark}

 \def\reg{\operatorname{reg}}
 \def\deg{\operatorname{deg}}
\def\adeg{\operatorname{adeg}}
\def\Gin{\operatorname{Gin}}
\def\gin{\operatorname{in}}

\def\mm{{\frak m}}
\def\pp{{\frak p}}

\begin{document}
\title{Finiteness of Hilbert functions and \\ bounds for
 Castelnuovo-Mumford regularity\\ of initial ideals}
\author{ L\^e Tu\^an Hoa }
\address{Institute of Mathematics\\ 18 Hoang Quoc Viet Road\\ 10307 Hanoi, Vietnam}
\email{lthoa@@math.ac.vn}
\thanks{The  author was  supported  in part by the National Basic Research
Program (Vietnam). The final preparation  of the article was done
during his stay at the Centre de Recerca Matematica (Spain).}
\keywords{Castelnuovo-Mumford regularity, local cohomology,
Hilbert function, Hilbert polynomial, initial ideal.}
 \subjclass{13D45, 13D40, 13P10}
\begin{abstract}  Bounds for the Castelnuovo-Mumford regularity and Hilbert coefficients
are given in terms of the arithmetic degree (if the ring is
reduced) or in terms of the defining degrees. From this it follows
that there exists only a finite number of Hilbert functions
associated to  reduced  algebras over an algebraically closed
field with a given arithmetic degree and dimension.  A good bound
is also given for the Castelnuovo-Mumford regularity of initial
ideals which  depends neither on  term orders nor on the
coordinates, and  holds for any  field.\end{abstract}
 \maketitle

\rightline{\it Dedicated to Professor J. Herzog}

\rightline{\it on the occassion of his 65-th birthday.}

\section*{Introduction}\smallskip

In the famous book SGA6, Kleiman  proved that given two positive
integers $e$ and $d$, there exists only a finite number of Hilbert
functions associated to  reduced and equidimensional $K$-algebras
$S$ over an algebraically closed field such that $\deg S  \le e$
and $\dim S =d$ (see [K, Corollary 6.11]).  An easier and
eleganter proof of this result can be found in a recent paper by
M. Rossi, N. V. Trung and G. Valla [RTV2]. Moreover, the paper
[RTV2] gives a rather general approach to derive the finiteness of
Hilbert functions. It is shown that this problem (for a certain
class of ideals) is equivalent to the boundness of the
Castelnuovo-Mumford regularity and the embedding dimension  (see
[RTV2, Theorem 2.3]). The first main purpose of this paper is to
extend Kleiman's result to  reduced $K$-algebras. A key point is
to find a suitable invariant to replace the degree. Of course, a
so-called extended degree is a choice, see [RVV, Corollary 4.4],
but such an invariant is very big.   It turns out that in our
situation one can take the so-called arithmetic degree - a notion
which maybe reflects better the complexity of ideals than the
usual degree (see [BM, Section 3] and [V, Chapter 9]).

 \begin{Theorem}\label{I1} Given two positive integers $a$ and $d$. Assume that $K$ is
 an algebraically closed field. Then there exists only a finite number of Hilbert
 functions associated to  reduced  $K$-algebras $S$  such that $\adeg S \le a$ and $\dim S =d$.
\end{Theorem}

Note that  the above result does not hold for an arbitrary algebra
(however see [RVV, Corollary 4.4] and [RTV2, Theorem 3.1] for a
possible generalization). As mentioned above, the main point in
the proof of Theorem \ref{I1} is to  bound the Castelnuovo-Mumford
regularity.  This is not hard to do (see Remark \ref{C6b}). However, a careful analysis allows us to  
establish the following explicit bound:

\begin{Theorem}\label{C4} Let $K$ be an arbitrary field and $I$  an arbitrary homogeneous ideal 
of $R= K[x_1,...,x_n]$. Assume that $S= R/I$ is a
reduced ring of dimension $d\ge 2$ and degree $e$. Then
$$\reg I \le (\frac{e(e-1)}{2} + \adeg I)^{2^{d-2}}.$$
\end{Theorem}

Applied to the case of reduced and equidimensional  algebras, the
bound of Theorem \ref{C4} is better than the one given in [RTV2,
Theorem 3.1 and Lemma 3.3]. In view of the Eisenbud-Goto
conjecture, the above bound is still too big, but it is a first
explicit bound stated in terms of the arithmetic degree. In order to
prove it, as in Lecture 14 of [M] (see also [K], [BM] and [RTV2]), we proceed by induction on the 
dimension. However, there is a different point: we simultaneously bound this invariant and the 
length of graded components of certain local cohomology modules  (see Theorem
\ref{C6} and also Theorem \ref{E5}). \par

The above technique can be also used to  estimate the
Castelnuovo-Mumford regularity of arbitrary homogeneous ideals in
terms of the maximal degree $\Delta$ of minimal generators of $I
\subset R=K[x_1,...,x_n]$. If $K$ is any field of zero
characteristic, from Giusti's paper [Gi] it follows that $\reg(I)
\le (2\Delta)^{2^{n-2}}$. Bayer and Mumford suggested that this
bound holds in any characteristic (see the comment after Theorem
3.7 in [BM]). Not long ago,  G. Caviglia and E. Sbarra proved that
this is indeed the case:
$$\reg I \le (\Delta^c + \Delta c - c+1 )^{2^{d-1}},$$
where $c= n-d$ (see [CS, Corollary 2.6]). In this paper we will give a completely
different proof for a slight improvement of this result (see Theorem \ref{E1}). \par

The next problem we are interested in is to give good bounds for
the Castelnuovo-Mumford regularity of initial ideals $\gin (I)$
with respect to {\it any term order} and in {\it any coordinates}.
Inspired by a result of Chardin and Moreno-Sosias, it was shown in
[HH] that if  $R/I$ is a Cohen-Macaulay ring of multiplicity $e\ge
2$ and $d \geq 2$, then $\reg(\gin(I)) \le
e^{2^{d-1}}/{2^{2^{d-2}}}$. Long before that M. Giusti [Gi] showed
that in characteristic zero we have $\reg \gin(I) \leq (2\Delta)^{2^{n-1} }$, provided the coordinates
are chosen generically and the term order is the lexicographic
order. Combining these facts with the above mentioned result of [CS] and [MM, II.2.2], one may ask 
whether  such a kind of bounds still
holds for any $\reg(\gin I)$. Our second main result confirms it:

\begin{Theorem}\label{I2} Let $K$ be an arbitrary field and and $I$  an arbitrary homogeneous 
ideal of $R$.
With respect to any term order and any coordinates we have
$$\reg(\gin I) \le (\frac{3}{2}\Delta^c+ \Delta)^{d2^{d-1}}.$$
Moreover, if $R/I$ is a reduced algebra, then we also have
$$\reg(\gin I) \le (\adeg I)^{(n-1)2^{d-1}}.$$
\end{Theorem}

  An immediate consequence of this theorem says that  the maximal degree
of a reduced Gr\"obner base, with respect to any term order and
any coordinates, is bounded by  $ (\frac{3}{2}\Delta^c+ \Delta)^{d2^{d-1}}$.  In
view of a remarkable example due to Mayr and Meyer, this bound is
nearly the best possible (see, e.g., Example 3.9 and Proposition
3.11 in [BM]).\par

In order to prove this theorem we develop further the method in
[HH]. Instead of initial ideals we consider a much bigger class:
the class of all ideals $J$ having the same Hilbert function as
$I$. So, although the title of the paper is about initial ideals, we are in fact dealing not much with 
them. However, by doing so one can use Gotzmann's regularity theorem to
bound $\reg J$ in terms of some data of $I$. Then, by virtue of
Theorem \ref{C4} and  Theorem \ref{E1}, we will see that
the only thing left is to estimate the Hilbert coefficients $e_i$
in terms of $\Delta$ or $\adeg (I)$ (see Lemmas \ref{F1} and
\ref{F3}). This problem is also of independent interest. Main
steps to do it may be explained as follows. First, using a recent
result by Herzog, Popescu and Vladoiu [HPV] one can bound
cohomological Hilbert functions (i.e. the length of graded
components of local cohomology modules) in terms of the
Castelnuovo-Mumford regularity. From that we  get  bounds for the
Hilbert coefficients by the Castelnuovo-Mumford regularity (see
theorems \ref{B3} and \ref{B3b}). The existence of such a bound
was predicted by [RTV2, Theorem 2.3], and this approach is
somewhat new, because usually one tries to estimate the latter
invariant by the former ones (see, e.g., [K] and [BrS, Section
17.2]). However, it is a surprising fact, that the relationships
between these invariants in theorems \ref{B3} and  \ref{B3b} are rather
simple. Let us give here a simple version of these results

\begin{Theorem}\label{I3} Let  $e_0=e,...,e_{d-1}$  be the Hilbert coefficents of $R/I$ and $b= 
\max\{ \Delta, \,  \adeg I\}$. Then
\begin{itemize}
\item[(i)] $|e_1|  \leq b^c\reg I$.
 \item[(ii)] For $i\ge 2$, $|e_i| \le \frac{3}{2} b^c ( \reg I)^i$.
\end{itemize}
\end{Theorem}

Combining theorems \ref{B3} and  \ref{B3b} with results on the Castelnuovo-Mumford
regularity found earlier we get bounds for $|e_i|$ in terms of
$\Delta$ or $\adeg (I)$ (see propositions \ref{C10} and \ref{E7}).
These bounds are huge: they are double exponential functions of
$i$. But they are good enough to prove
Theorem \ref{I2}. Furthermore, theorems \ref{B3} and \ref{B3b}
sometimes give really good bounds for
 $|e_i|$ if we already know a good estimation for the Castelnuovo-Mumford regularity
  (see corollaries \ref{B6} and \ref{B7}).\par

We now give a brief content of the paper. We
 prove Theorem \ref{C4} in Section \ref{C}, and reprove in Section \ref{E} the
Caviglia-Sbarra bound on the Castelnuovo-Mumford regularity of an
arbitrary homogeneous ideal in terms of the degrees of its
defining equations (Theorem \ref{E1}). Section \ref{A} is devoted
to bounding Hilbert cohomological functions in terms of the
Castelnuovo-Mumford regularity (see Theorem \ref{A4}). Bounds on
Hilbert coefficients are given in Section \ref{B}.   Putting
results of the sections \ref{C} and \ref{B} together we are able
to prove Theorem \ref{I1} without using [RTV2]. This is done
in Section \ref{D}.  Theorem \ref{I2} is proved in the last
Section \ref{F}. We refer the readers to Eisenbud's book [E] for
unexplained terminology. \vskip0.5cm

\noindent {\bf Acknowledgment}:  The author would like to thank
the Centre de Recerca Matematica (Spain) for the financial support and
hospitality during the final preparation of this article. He is grateful to the referee
for his/her useful remarks and suggestions which lead to an improvement of some main results of the 
paper.

\section{Bounds in terms of the arithmetic degree} \label{C}\smallskip

Throughout this paper, if not otherwise stated, $K$ is an
arbitrary  field, $R= K[x_1,...,x_n]$ is a polynomial ring
and $I \subset R$ is a homogeneous ideal of dimension $d$. However, all invariants considered in 
this paper are not changed under passing from $K$ to $K(u)$, where $u$ is  a new indeterminates. 
Hence, in  proofs we may always assume that $K$ is an infinite field. This assumption guarantees 
the choice of generic elements.

 Let $c=
n-d$. Note that $c$ is the true codimension of $I$ if $I$ does not
contain a linear form. Let $\mm = (x_1,...,x_n)$ denote the
maximal homogeneous ideal of $R$ and set $S= R/I$. Let us recall
some notions.\par

For an artinian $\Bbb Z$-graded  module $N$, let
$$\text{\rm end}(N)  = \max \{ t;\ N_t \neq 0\}$$
(with the convention $\max \emptyset = - \infty$). Further, let
$$a_i(R/I) = \text{\rm end}( H^i_\mm(R/I)),$$
where $H^i_\mm(R/I)$ is the local cohomology module with the
support in $\mm$ and  $ 0 \leq i\leq d$.
 The {\it Castelnuovo-Mumford regularity} is the number $$\reg(R/I) =
 \max\{ a_i(R/I) + i;\ 0 \leq i\leq d\}.$$
Note that $\reg(I) = \reg(R/I) +1$. Sometimes we also use the
notation
\begin{eqnarray} \label{reg1}
\reg_k(R/I) =
 \max\{ a_i(R/I) + i;\ k \leq i\leq d\},
\end{eqnarray}
where $k$ is a non-negative integer.\par

Following Brodmann and Sharp [BrS], the function
$$h^i_S(t) := \ell (H^i_\mm(S)_t)$$
is called the $i$-th {\it Hilbert homological function} of $S$,
where $\ell(.)$ denotes the dimension of a vector space over $K$.
Let  $H_S(t)$ and $P_S(t)$ denote the Hilbert function and the
Hilbert polynomial of $S$, respectively. We will often use the
Grothendieck-Serre  formula
\begin{eqnarray} \label{EB1}
P_S(t) - H_S(t) = \sum_{i=0}^d (-1)^{i+1}h^i_S(t).
\end{eqnarray}\par

The leading coefficient of $P_S(t)$, multiplied by $(d-1)!$, is called the degree of $S$ and denoted 
by $\deg S$. We also denote $\deg S$ by $e(S)$, or just by $e$.
The arithmetic degree  is defined as follows:
$$\adeg S = \adeg I = \sum_{\pp \in \text{Ass}(R/I)} \ell(H^0_{\mm_\pp}(R_\pp/I_\pp)) 
\deg(R/\pp),$$
(see [BM, Definition 3.4] and [V, Definition 9.13]). The number $
\ell(H^0_{\mm_\pp}(R_\pp/I_\pp))$ is the multiplicity of the component
$\pp$ with respect to $I$. In this definition $\pp$ runs over all
associated primes of $S$, while the usual degree $\deg S$ can be
computed by a similar formula, but the sum is only taken over
primes of the highest dimension. Thus
$$\adeg S \ge \deg S,$$
and the equality holds if and only if $S$ is a pure-dimensional
ring.  \par

In this section we prove Theorem \ref{C4}. We need some auxiliary
results.

\begin{Lemma}\label{C1} Let $S$ be an one-dimensional Cohen-Macaulay ring.
Then
$$h^1_S(0) + \cdots + h^1_S(\reg S-1) \le e(e-1)/2.$$
\end{Lemma}

\begin{pf} Since $P_S(t) = e$, from the Grothendieck-Serre formula (\ref{EB1}) we have
$$h^1_S(t) = e - H_S(t).$$
 Let $r= \reg S$. Since $S$ is a Cohen-Macaulay ring, its Hilbert-Poincare series can be
 written in the  form
$$HP_{S}(z) := \sum_{i\ge 0} H_{S}(i)z^i = \frac{1 + h_1z + \cdots + h_rz^r}{1-z},$$
where $h_1,...,h_r$ are positive integers (see, e.g., [V, p.
240]). From this it follows that
$$H_S(t) = 1+ h_1 + \cdots + h_t \ge t+1$$
for all $t\le r$.
 Moreover, under the Cohen-Macaulay assumption,  $r \le e-1$ . Hence
$$h^1_S(0) + \cdots + h^1_S(\reg S-1) \le re - (1+ \cdots + r) = r(2e-r-1)/2 \le e(e-1)/2.$$
\end{pf}

\begin{Lemma}\label{C2}  Assume that $S= R/I$ is a reduced ring of dimension at least two. Then
$$h^1_S(-1) \le \adeg I - e.$$
\end{Lemma}

\begin{pf} Since $S$ is reduced, one may write $I = J \cap Q$, where $J$ is the intersection
of all associated primes of $R/I$ of dimension at least 2, and $Q$
is the intersection of all associated primes of $R/I$ of dimension
1. By [HSV, Lemma 1] we have $h^1_{R/J}(-1) = 0$. Thus if $Q=R$,
then $h^1_S(-1) =0$. Assume that $Q \neq R$. Since $J \neq R$ and
$R/I$ has no embedded primes, $J+Q $ is an $\mm$-primary ideal,
i.e. $\dim R/(J+Q) = 0$. The exact sequence
$$0 \to S \to R/J \oplus R/Q \to R/(J+Q) \to 0$$
implies
$$h^1_S(-1) = h^1_{R/J}(-1) +  h^1_{R/Q}(-1) =  h^1_{R/Q}(-1).$$
Note that $\deg R/Q = \adeg I - \adeg J \le \adeg I -e.$ Since
$R/Q$ is an one-dimensional ring, by the Grothendieck-Serre
formula, we have
$$h^1_{R/Q}(-1) = \deg R/Q \le \adeg I -e.$$
\end{pf}

The proof of Theorem \ref{C4} is proceeded by induction. The next
two lemmas allow us to do induction. The first one is concerning
the behavior of the arithmetic degree by hyperplane section. It is
more subtle than the usual degree, see [MVY]. However we have

\begin{Lemma}\label{B2} Let $K$ be an infinite field, and $S=R/I$  an arbitrary  ring of 
dimension at least
two and positive depth. Assume  that $x_n$ is chosen generically.
Let $T=R/((I,x_n):\mm^{\infty})$ and $r= \reg T$. Then:
\begin{itemize}
\item[(i)] $\reg T \le \reg S$.
\item[(ii)] $\adeg T \le \adeg S$.
\end{itemize}
\end{Lemma}

\begin{pf}
  (i) Since $x_n$ is generic, it is a regular element on $S$. We have
$$\reg T  = \reg_1 S/x_nS \le \reg S/x_nS =  \reg S.$$

(ii) For an $R$-module $M$ and $r\ge -1$, let
$$\adeg_r (M) = \sum_{\pp \in \text{Ass}(M),\ \dim R/\pp = r+1} \ell(H^0_{\mm_\pp}(M_\pp))
\deg(R/\pp)$$ (see [BM, Definition 3.4]). Since $x_n$ is generic, by
the prime avoidance lemma,
 we may  assume
$$x_n \not\in \cup\{ \pp;\ \mm \neq \pp \in \text{Ass}(S)\cup_{j\ge 1}
 \text{Ass}( \text{Ext}^{n-j}_{R}(S,R))\}.$$
  By [MVY, Corollary 2.5] it
follows that
$$\adeg_{r-1}(T) = \adeg_r(S)  \ \ \text{for\ all} \ r \ge 1.$$
Since $S$ and $T$ have no zero-dimensional component, we get
$$\begin{array}{ll}
\adeg T  &= \adeg_0(T) + \cdots + \adeg_{d-1}(T)\\
&=  \adeg_1(S) +\cdots +  \adeg_d(S) \le  \adeg S.
 \end{array}$$
\end{pf}

 The first three
statements of the next lemma are contained in the proof of Mumford's theorem on page 101 of the 
book [M] 
(cf. also [K, Proposition 1.4], [RTV1, Theorem 1.4] and [RTV2, Theorem 1.3]). In order
to make the paper more self-contained, we give here a sketch of
the proof. The proof of (iii) here is also simpler.

\begin{Lemma}\label{C5} Let $K$ be an infinite field and $S= R/I$  a reduced ring of dimension 
at least two.
Assume that $x_n$ is chosen generically. Let
$T=R/((I,x_n):\mm^{\infty})$ and $r= \reg T$. Then $T$ is also a
reduced ring and we have
\begin{itemize}
\item[(i)] $\reg_2(S) \le r$ (see the definition in (\ref{reg1})).
\item[(ii)] $h^1_S(t) \ge h^1_S(t+1)$ for all $t \ge r-1$.
\item[(iii)] $\reg S \le r + h^1_S(r-1)$. \item[(iv)] $h^1_S(t)
\le h^1_T(0) + \cdots + h^1_T(t) + \adeg I - e$, for all $t\ge 0$.
\end{itemize}
\end{Lemma}

\begin{pf} Note that $T$ can be considered as the homogeneous coordinate ring of a generic 
hyperplane
section of the scheme $\text{Proj}(S)$. Since $K$ is an infinite
ring and $x_{n}$ is generic, by    Bertini's theorem [FOV,
Corollary 3.4.14] it follows that $T$ is reduced.

The long exact sequence
\begin{eqnarray}
0 \to   H^0_{\mm}(S/x_nS)_{t} & \to &  H^1_{\mm}(S)_{t-1} \to  H^1_{\mm}(S)_{t}
{\overset{\varphi_t}{\longrightarrow}}  H^1_{\mm}(S/x_nS)_t  =  H^1_{\mm}(T)_{t} 
\label{EC1}\\
 & \to & H^2_{\mm}(S)_{t-1}\to H^2_{\mm}(S)_{t} \to  \cdots \nonumber
\end{eqnarray}
implies (i) and the short exact sequence
$$0 \to   H^0_{\mm}(S/x_nS)_{t} \to  H^1_{\mm}(S)_{t-1} \to  H^1_{\mm}(S)_{t} \to 0$$
for all $t \ge r$. This yields (ii). If $h^1_S(t_0-1) \ge
h^1_S(t_0)$ for some $t_0 \ge r+1$,  we would have
$h^0_{S/x_nS}(t_0) = 0$. Since $\reg_1(S/x_nS) = \reg T = r$, it
then implies that $\reg (S/x_nS) \le t_0$.  Hence $h^1_S(t_0) =
h^1_S(t_0+1) = \cdots = 0$. Therefore $h^1_S(t)$ is strictly
decreasing to zero when $t \ge r$, which implies (iii).

 It remains to show (iv). From the exact sequence (\ref{EC1}) we have
$$h^1_S(u) - h^1_S(u-1) = \ell( \text{Im}(\varphi_u)) - h^0_{S/x_nS}(u) \le h^1_T(u)$$
for all $u \in \Bbb Z$. Adding these inequalities and using
Lemma \ref{C2} we get
$$\begin{array}{ll}
h^1_S(t) & \le h^1_T(0) + \cdots + h^1_T(t) + h^1_S(-1) \\
& \le h^1_T(0) + \cdots + h^1_T(t) + \adeg I - e
 \end{array}$$
for all $t\ge 0$.
\end{pf}

Theorem \ref{C4} is a part of the following result. If $ a\in \Bbb R$, we denote by $[a]$
the largest integer not exceeding $a$.

\begin{Theorem}\label{C6} Assume that $S= R/I$ is a reduced ring of dimension at least two. Let
$$m= \frac{e(e-1)}{2} + \adeg I.$$
Then
\begin{itemize}
\item[(i)] $\reg S \le m^{2^{d-2}} -1$. \item[(ii)] For all $t\ge
0$, we have $h^1_S(t) \le m^{2^{d-2}} - e\cdot m^{[2^{d-3}]} $.
\end{itemize}
\end{Theorem}

\begin{pf} We may assume that $x_n$ is generic and choose $T$ as in the previous lemma.
Hence $T$ is a reduced ring. Set  $r= \reg T$.

Let $d=2$. In order to show (ii), by Lemma \ref{C5}(ii), we may
assume that $t\le r-1$. Note that $T$ is a Cohen-Macaulay ring and
$e(T) = e$. Then Lemma \ref{C5}(iv) and Lemma \ref{C1} yield:
$$\begin{array}{ll}
h^1_S(t) & \le h^1_T(0) + \cdots + h^1_T(t) + \adeg I - e \\
& \le h^1_T(0) + \cdots + h^1_T(r-1) + \adeg I - e \\
& \le  \frac{e(e-1)}{2} + \adeg I - e = m - e.
 \end{array}$$
Using this inequality and the fact that $r \le e-1$ (since $T$ is
a Cohen-Macaulay ring), by Lemma \ref{C5}(iii) we get
$$\reg S \le e-1 + m-e = m-1.$$
Thus the case $d=2$ is proven.

Let $d\ge 3$. Since $\dim T = d-1, \ e(T) = e$, and $\adeg T \le
\adeg S$ (by Lemma \ref{B2}(ii)), the induction hypothesis gives
\begin{eqnarray}\label{EC2}
r \le m^{2^{d-3}} -1,
\end{eqnarray}
and for all $t\ge 0$
\begin{eqnarray}\label{EC3}
h^1_T(t) \le m^{2^{d-3}} - e\cdot m^{[2^{d-4}]} \le m^{2^{d-3}} -e
.
\end{eqnarray}
In order to prove (ii), again by Lemma \ref{C5}(ii), we may assume
that $t\le r-1$. Then, by Lemma \ref{C5}(iv), for all $t\ge 0$ we
have
$$\begin{array}{ll}
h^1_S(t) & \le h^1_T(0) + \cdots + h^1_T(t) + \adeg I - e \\
& \le r(m^{2^{d-3}} -e)+ \adeg I - e  \ \ \text{(by \  (\ref{EC3}))} \\
& \le  (m^{2^{d-3}} -1)(m^{2^{d-3}} -e) + \adeg I - e \ \ \text{(by \  (\ref{EC2}))} \\
& = m^{2^{d-2}} - e\cdot m^{2^{d-3}} - m^{2^{d-3}} + \adeg I\\
& \le m^{2^{d-2}} - e\cdot m^{2^{d-3}}.
 \end{array}$$
To prove (i) we use (ii) and Lemma \ref{C5}(iii) :
$$\reg S \le r + h^1_S(r-1) \le  m^{2^{d-3}} -1 + m^{2^{d-2}} - e\cdot m^{2^{d-3}} \le 
m^{2^{d-2}}-1.$$
\end{pf}

\begin{Remark}\label{C6b} {\rm  In order to get a bound for $\reg S$ in terms of $\adeg S$ and 
$d$, by induction on $d$, it suffices to estimate $h^1_S(t)$ for $t\geq 0$. This can be  easily done 
by using the following well-known inequality
$$h^1_S(t) + H_S(t) \leq \adeg S {t+d-1 \choose d-1}$$
(recall that $S$ is reduced). This was pointed out by the referee. However a direct application of this 
inequality would only lead to a bound of the following type 
$$\reg I \leq (\adeg I)^{d-1)!}.$$}
\end{Remark}

If $K$ is an algebraically closed field, then a result of Gruson,
Lazarsfeld and Peskine  [GLP] yields a better bound for the case
$d=2$ as shown in the following statement.

\begin{Proposition}\label{C7} Let $K$ be an algebraically closed field. Assume that $R/I$
is a reduced ring of dimension two. Then $\reg I \le \adeg I$.
\end{Proposition}

\begin{pf} Write $I = J \cap Q$ as in the proof of Lemma \ref{C2}. Then we have an exact sequence:
$$0 \to H^0_{\mm}(R/J+Q)_t \to H^1_{\mm}(R/I)_t \to H^1_{\mm}(R/J)_t \oplus 
H^1_{\mm}(R/Q)_t \to 0,$$
and
$$H^2_{\mm}(R/I)_t \cong H^2_{\mm}(R/J)_t \oplus H^2_{\mm}(R/Q)_t.$$
By [GLP, Theorem 1.1] (see also Remark on p. 497 there), $\reg R/J
\le e-1 \le \adeg I -1$. So we may assume that $Q \neq R$. Since
$R/Q$ is one-dimensional and reduced, it is a Cohen-Macaulay ring.
Hence $\reg R/Q \le \adeg R/Q -1 < \adeg I$. To complete the proof
it suffices to show that
$$H^0_{\mm}(R/J+Q)_t = 0 \ \ \text{for \ all}\ t \ge \adeg I -1.$$
Since $\reg J \le e$, $J$ is generated by elements of degree $\le
e$. Hence one may choose an element $x\in J$ of degree $e$ such
that $x$ does not belong to any prime in $Q$, i.e. $x$ is a
regular element on $R/Q$. Then we have
$$\reg R/(Q,x) = \reg R/Q + e-1 \le \deg R/Q -1 + e-1 = \adeg I -2.$$
Since $R/(Q,x)$ is a zero-dimensional ring, this means
$(R/(Q,x))_t =0$ for all $t \ge \adeg I -1$. The inequality
$\ell((R/J+Q)_t) \le \ell((R/(Q,x))_t)$ gives us $H^0_{\mm}(R/J+Q)_t =
(R/J+Q)_t  = 0$ for all $t \ge \adeg I -1$, as required.
\end{pf}

The above proposition says that in dimension two one can replace $m$ by $\adeg S$ in Theorem 
\ref{C6}(i). However this does not work for Theorem \ref{C6}(ii) as shown by the following example.

\begin{Example}\label{C8} {\rm Given $e\ge 6$ and $S= K[x^e, x^{e-1}y, xy^{e-1}, y^e]$.
Then for $0\le t \le e-2$ one can show that $h^1_S(t) = te+1 - (t+1)^2$, while $\adeg S - e = 0$. 
Taking $t_0 = [\frac{e-2}{2}]$, one can see that  $h^1_S(t_0)$ is approximately a half of the bound 
in Theorem \ref{C6}(ii).
}\end{Example}

Theorem  \ref{C4} does not hold if the ring $R/I$ is not reduced.

\begin{Example}\label{C9} {\rm (see [V, Example 9.3.1]) Let $S= K[x,y,u,v]/((x,y)^2,
xu^t+yv^t),\ t\ge 1$. Then  $\adeg S = e = 2$, while $\reg S = t$ can be arbitrarily large.
}\end{Example}

 \section{Bounds in terms of degrees of defining equations} \label{E}\smallskip

In this section we study arbitrary homogeneous ideals. We will
always write the degrees of polynomials in a minimal homogeneous
basis of $I$  in a decreasing sequence
$$\Delta := \delta_1 \ge  \delta_2 \ge \cdots $$
 and assume $\Delta \ge 2$. As mentioned in the introduction,
 G. Caviglia and E. Sbarra already proved that
$$\reg I \le (\Delta^c + \Delta c - c+1 )^{2^{d-1}}$$
(see [CS, Corollary 2.6]). The purpose of this section is to prove
the following theorem which is a slight improvement of the above
result.

\begin{Theorem}\label{E1} Let $K$ be an arbitrary field and $I$ be an arbitrary homogeneous
 ideal of  dimension  $d\ge 1$. Then
$$\reg I \le (\delta_1 \cdots  \delta_c + \Delta -1)^{2^{d-1}} 
\le (\Delta^c +\Delta -1)^{2^{d-1}}.$$
\end{Theorem}

The proof of [CS] uses properties of Borel-fixed ideals. The proof
here is completely different and simpler than the one in [CS]. The main
idea of the proof  is similar to that of Theorem \ref{C4}.  We need some technical lemmas. For
short, set
$$\sigma =  \delta_1 +\cdots+  \delta_c - c \ \ \text{and} \ \ \pi =
\delta_1 \cdots  \delta_c .$$

If $S= R/I$, then  we also write $\Delta = \Delta(S),\
 \delta_1 =  \delta_1(S), ...$ to emphasize their dependence on $S$ (or $I$). The 
following result was pointed out by the referee to the author.  Subsequently,  it slightly improves our 
original Theorem \ref{E5} and Theorem \ref{F4}

\begin{Lemma} {\rm [Sj, Theorem 2]}\label{E2}   If $\dim S \le 1$, then $\reg S \le \sigma + \Delta 
-1$.
\end{Lemma}

The next result is a special case of [HH, Lemma 3].

\begin{Lemma}\label{E3} Assume that $d=1$. Then for all $t\ge 1$, $h^0_S(t) \le \pi -1$.
\end{Lemma}

Recall that an element $x\in \mm$ is called {\it filter regular}
if $0: \mm^\infty$ is of finite length.

\begin{Lemma}\label{E4} Assume that $\dim S\ge 1$ and $x=x_n$ is a filter regular
element on $S$. Let $T=S/xS$ and $r \ge \max\{ \reg T,\ \Delta-1\}$.  Then
\begin{itemize}
\item[(i)] $\reg_1(S) \le r$ (see the definition in (\ref{reg1})).
\item[(ii)] $h^0_S(t) \ge h^0_S(t+1)$ for all $t \ge r$.
\item[(iii)] $\reg S \le r + h^0_S(r)$. \item[(iv)] $h^0_S(t) \le
h^0_T(1) + \cdots + h^0_T(t) $, for all $t\ge 1$.
\end{itemize}
\end{Lemma}

\begin{pf} (i)-(iii) were shown in the proof of [BM, Proposition 3.8]. It follows from
the following exact sequence
$$ 0 \to   (0:x)_{t-1} \to H^0_{\mm}(S)_{t-1} \to H^0_{\mm}(S)_{t}
{\overset{\varphi_t}{\longrightarrow}} H^0_{\mm}(T)_{t} \to  H^1_{\mm}(S)_{t-1}
 \to  H^1_{\mm}(S)_{t} \to    \cdots $$
For (iii) we need also the assumption $r\ge \Delta -1$ in order to
apply the regularity criterion of [BS, Theorem 1.10].

From the above exact sequence we have
$$h^0_S(u) - h^0_S(u-1) = \ell( \text{Im}(\varphi_u)) - \ell((0:x)_{u-1}) \le h^0_T(u)$$
for all $u \in \Bbb Z$. Since $h^0_S(0) = 0$, adding these
inequalities gives us (iv).
\end{pf}

Theorem \ref{E1} is a part of the following

\begin{Theorem}\label{E5} Let $d\ge 1$. Then
\begin{itemize}
\item[(i)] $\reg S \le  (\pi + \Delta -1)^{2^{d-1}} -1$. 
\item[(ii)]
For all $t\ge 1$, we have $h^0_S(t) \le  (\pi + \Delta -1)^{2^{d-1}} -
 (\pi + \Delta -1)^{[2^{d-2}]} $.
\end{itemize}
\end{Theorem}

\begin{pf} Keep the notation of Lemma \ref{E4}. Let $I'$ denote the image of $I$
in $K[x_1,...,x_n]/(x_n) \cong K[x_1,...,x_{n-1}] =: R'$. Then $T \cong R'/I'$
and it is clear that $\sigma(I') \le \sigma$ and $\pi(I') \le \pi$.

First let $d=1$. With the above remark, Lemma \ref{E3} yields
$h^0_S(t) \le \pi -1$ for all $t\ge 1$. Thus (ii) holds. By Lemma
\ref{E2}, 
$$\reg S \le \delta_1 + \delta_2 +\cdots + \delta_c - c +\Delta -1.$$
By induction on $c$ we  get 
$$\reg S \le (\delta_1\cdots \delta_{c-1} -1) + (\delta_c -1) + \Delta -1 \le \delta_1\cdots \delta_c -1 
+ \Delta -1 = \pi +\Delta - 2.$$
 Thus the case $d=1$ is proven.

Now let $d\ge 2$. With the remark at the beginning of the proof
and by the induction hypothesis we may assume that
$$ \reg T \le  (\pi + \Delta -1)^{2^{d-2}} -1,$$ and
$$h^0_T(t) \le  (\pi + \Delta -1)^{2^{d-2}} -  (\pi + \Delta -1)^{[2^{d-3}]}
\le  (\pi + \Delta -1)^{2^{d-2}} -1. $$
Let $r= (\pi + \Delta -1)^{2^{d-2}} $. Obviously $r \ge \Delta $. In
order to prove (i), by Lemma \ref{E4}(ii), we may assume that
$t\le r$. Then,  by Lemma \ref{E4} (iv) and the induction
hypothesis, we have
$$\begin{array}{ll} h^0_S(t) & \le r ((\pi + \Delta -1)^{2^{d-2}} -1) \le  
(\pi + \Delta -1)^{2^{d-2}} ( (\pi + \Delta -1)^{2^{d-2}} -1) \\
&=  (\pi + \Delta -1)^{2^{d-1}} - (\pi + \Delta -1)^{2^{d-2}}. 
\end{array}$$
Thus (ii) is proven. Using this and Lemma \ref{E4}(iii) we
immediately get (i).
 \end{pf}

\begin{Remark}\label{E6} {\rm  a) The bound in Theorem \ref{E1} is nearly the best possible.
It was shown  that there is an ideal $I$, due to Mayr and Meyer, generated by $10n-6$
forms of degree at most 4 in $10n +1$ variables such that
$\reg(I) > 4^{2^{n-1}} + 1$ (see, e.g.,    [BM, Example 3.9 and
Proposition 3.11]). \par

b) In this paper we are not interested in giving the best possible bounds for $\reg S$ which are then 
more complicated to formulate. On the other hand, for rings of small dimension, there are 
also some bounds which are much better than the ones in Theorem \ref{E1}.  See, e.g., a recent paper 
[CF] for $d\leq 2$. However an application of such results  to our proof does not significantly improve 
the bound in Theorem \ref{E5} for a larger 
$d$.}\end{Remark}

\section{Hilbert cohomological functions} \label{A}\smallskip

In this section we give a bound on $h^i_S(t)$.  First, we do this
for Borel-fixed ideals. We need some notation and results from
[HPV]. Let $I\neq 0$ be a monomial ideal. Denote by $G(I)$ the
unique set of monomial generators of $I$. For a monomial $u$, let
$m(u)$ be the maximal index of a variable appeared in $u$. Set
$$m(I) = \max \{ m(u); \ u \in G(I)\}.$$
We recursively define an ascending chain of monomial ideals
$$I = I_0 \subset I_1 \subset \cdots \subset I_{l+1} = R$$
as follows: let $I_0=I$. Suppose $I_j$ is already defined. If $I_j
= R$, then the chain ends. Otherwise, let $n_j = m(I_j)$ and set
$$I_{j+1} = I_j : x_{n_j}^\infty := \cup_{k=1}^\infty I_j : x_{n_j}^k.$$

A stable ideal under the action of upper triangle matrices is
called {\it Borel-fixed}. It is always a monomial ideal. If $I$ is
a Borel-fixed ideal, then $(x_1,...,x_c)$ is the unique minimal
associated prime of $R/I$ (see [E, Corollary 15.25]). Hence in
this case $n \geq n_0 > n_1 > \cdots > n_l =c$. For $j=0,...,l$,
let $J_j \subset K[x_1,...,x_{n_j}]$ be the monomial ideal with
$G(I_j)=G(J_j)$. Denote by
$$J_j^{sat} = J_j :
(x_1,...,x_{n_j})^\infty$$
 the {\it saturation} of $J_j$. Then by [HPV,
Corollary 2.6] and local duality we have

\begin{Lemma}\label{A1} Let $I\neq 0$ be a Borel-fixed ideal. Then $H^j_\mm(S) = 0$
 if $j \not\in \{n-n_0,...,n-n_l\}$, and we have an isomorphism of $\Bbb Z$-graded
 $R$-modules: $$H^j_\mm(S) \cong (J_i^{sat}/J_i)[x_{n_i+1}^{-1},...,x_n^{-1}],$$
if $j=n-n_i$ for some $i=0,...,l$.
\end{Lemma}

In the sequel, for a Borel-fixed ideal $I$ let us denote
\begin{eqnarray}\label{EA1}
B := B(I) = \ell(R/(I,x_{c+1},...,x_n)).
\end{eqnarray}
For short, set $e = \deg(I)$. Note that $B \geq e$.

\begin{Lemma}\label{A2} Let $I\neq 0$ be a Borel-fixed ideal. Then
\begin{itemize}
\item[(i)] $\ell(J_l^{sat}/J_l) =e$. \item[(ii)] For $i< l$ and
all $t \ge 0$ we have
$$\ell([J_i^{sat}/J_i]_t) \le (B-1){t+n_i-c-2\choose n_i-c-1}.$$
\end{itemize}
\end{Lemma}

\begin{pf}  Let $M= J_i^{sat}/J_i$ and  $R' = K[x_1,...,x_c]$. Since $(x_1,...,x_c)$
is the unique minimal associated prime of $R/I$, by the
construction we have $I \subseteq I_i \subseteq (x_1,...,x_c)$.

 Let $i=l$. We have $J_l^{sat} = R'$ and $J_l = I: (x_{c+1},...,x_{n})^\infty$. Hence
$$\ell(M) = \ell(R'/I: (x_{c+1},...,x_{n})^\infty) = \ell((R/I)_{(x_{c+1},...,x_{n})}) = e.$$

Let $i<l$. Set $R" =K[x_1,...,x_{n_i}]$. By the definition $J_i =
G(I_i)R"$. Hence $x_{c+1},...,x_{n_i}$ is a s.o.p. of $R"/J_i$,
and $ I \cap R' \subseteq J_i \cap R'.$ This implies
\begin{eqnarray}\label{EA0}
\ell(\frac{R"}{(J_i, x_{c+1},...,x_{n_i})}) = \ell(\frac{R'}{J_i
\cap R'}) \le \ell(\frac{R'}{I\cap R'}) = \ell(\frac{R}{(I,
x_{c+1},...,x_{n})}) = B.
\end{eqnarray}
On the other side, the inclusion $J_i^{sat} \subseteq
(x_1,...,x_c)$ yields
$$\ell(M_t) \leq \ell(\frac{(x_1,...,x_c)R"}{J_i}) = \ell(\frac{R"}{J_i})
 - \ell(\frac{R"}{(x_1,...,x_c)R"} =
\ell(\frac{R"}{J_i}) - {t+n_i-c-1\choose n_i-c-1}.$$ By [RVV,
Proposition 2.4] and (\ref{EA0}) we have
$$\begin{array}{ll}
\ell(\frac{R"}{J_i}) &\le (\ell( \frac{R"}{(J_i, x_{c+1},...,x_{n_i})}) -1)
{t+n_i-c-2\choose n_i-c-1} + {t+n_i-c-1\choose n_i-c-1}\\
&\le (B-1){t+n_i-c-2\choose n_i-c-1} + {t+n_i-c-1\choose n_i-c-1}.
 \end{array}$$
Hence $\ell(M_t) \le (B-1){t+n_i-c-2\choose n_i-c-1}$.
\end{pf}

\begin{Lemma}\label{A3} Let $I\neq 0$ be a Borel-fixed ideal and $S=R/I$. Then
\begin{itemize}
\item[(i)] $h^0_S(t) \le (B-1){t+d-2\choose d-1} $ for all $t\ge
0$. \item[(ii)] For $1 \le j \le d-1$ and $t \le \reg S$:
$$h^j_S(t) \le (B-1){\reg S+d-j-2\choose d-j-1} {\reg S-t\choose j}.$$
\item[(iii)] For $t< \reg S$:
$$h^d_S(t) \le e{\reg S-t-1\choose d-1}  \le B{\reg S-t-1\choose d-1}.$$
\end{itemize}
\end{Lemma}

\begin{pf} By virtue of Lemma \ref{A1}, (i) is a special case of Lemma \ref{A2}
(when $i=0$). Let $j\ge 1$. By Lemma \ref{A1} we may assume that
$j=n-n_i$ for some $i>0$. Let $M = J_i^{sat}/J_i$. Lemma \ref{A1}
implies
\begin{eqnarray}
h^j_S(t) &=& \sum_{u= 1}^{\text{\rm end}(M)} \ell(M_u){u-j-t+j-1\choose j-1}  \label{EA2}\\
&\le& [\max_{1\le u\le \text{\rm end}(M)}\ell(M_u)]  \sum_{v=
1}^{\text{\rm end}(M)-j-t}{v+j-1\choose j} \nonumber\\&=&
[\max_{1\le u\le \text{\rm end}(M)}\ell(M_u)]  {\text{\rm
end}(M)-t\choose j}. \label{EA3}
\end{eqnarray}
(In the above calculation we set ${a\choose b} =0$ if $b\geq 0$
and $a<b$.) Moreover, again by Lemma \ref{A1}, $\text{\rm end}(M)
= a_j(S) + j \le \reg S$ (see also [HPV, Corollary 2.7]). Since
$n_i-c = n-j-c = d-j$, Lemma \ref{A2}(ii) yields
$$\max_{1\le u\le \text{\rm end}(M)}\ell(M_u) \le (B-1) \max_{1\le u\le \reg S}
{u+d-j-2\choose d-j-1} =
(B-1) {\reg S+d-j-2\choose d-j-1}.$$ From this and (\ref{EA3}) we
get (ii).

Let $j=d$. Then $n_i=c$. From (\ref{EA2}) we have
$$\begin{array}{ll}
h^d_S(t) & \le [\max_{1\le u\le \text{\rm end}(M)}{u - t -1\choose d-1}]
 \cdot  \sum_{u= 1}^{\text{\rm end}(M)} \ell(M_u)\\
& \le {\reg S-t-1\choose d-1} \cdot \ell(M)\\
&= e{\reg S-t-1\choose d-1} \ \ \text{\rm (by \ Lemma \
\ref{A2}(i))}.
 \end{array}$$
\end{pf}

Now we can bound the Hilbert cohomological functions of an
arbitrary homogeneous ideal $I$. Recall that the defining degrees
of $I$ are written in a decreasing sequence
$$\Delta := \delta_1 \ge  \delta_2 \ge \cdots,$$
and assume $\Delta \ge 2$.

In the proof of the following theorem, a result of Vasconcelos on
the reduction number plays an essential role. We use initial ideals in order to go back to the 
situation of the previous result.

\begin{Theorem}\label{A4} Let $I$ be an arbitrary homogeneous ideal of $R$ and $S=R/I$. Let
$$b= \min\{  \delta_1 \cdots  \delta_c,\ (\adeg I)^c\}.$$
Then
\begin{itemize}
\item[(i)] $h^0_S(t) \le (b-1){t+d-2\choose d-1} $ for all $t\ge
0$. \item[(ii)] For $1 \le j \le d-1$ and $t \le \reg S$:
$$h^j_S(t) \le (b-1){\reg S+d-j-2\choose d-j-1} {\reg S-t\choose j}.$$
\item[(iii)] For $t< \reg S$:
$$h^d_S(t) \le e{\reg S-t-1\choose d-1}  \le b{\reg S- t -1\choose d-1}.$$
\end{itemize}
\end{Theorem}

\begin{pf} Let $\Gin I$ denote the generic initial ideal of $I$ with respect to
 the reverse lexicographic order. Then $\Gin I$ is a Borel-fixed ideal. Moreover we
 may assume that the coordinates  $x_1,..,x_n$ are chosen generically. By
  [BS, Lemma 2.2 and Theorem 2.4] we have
$$\ell(R/(I, x_{c+1},...,x_{n})) = \ell(R/(\Gin I, x_{c+1},...,x_{n})),$$
and
$$\reg(R/I) = \reg(R/\Gin I).$$
By Macaulay's theorem: $e(R/I) = e(R/\Gin I)$. Moreover, by [S,
Theorem 2.4]
$$h^i_{R/I}(t) \le h^i_{R/\Gin I}(t) $$
for all $i\ge 0$ and $t\in \Bbb Z$. Hence, the theorem immediately
follows from the previous lemma if we can show that
$$ B:= \ell(R/(I, x_{c+1},...,x_{n})) \le b.$$

a) Let $I'$ denote the image of $I$ in $R' = R/(x_{c+1},...,x_{n})
\cong K[x_1,...,x_c]$.  It is a $(x_1,...,x_c)$-primary ideal.
Since $I$ can be generated by elements of degrees  $d_1\le
\delta_1, d_2\le \delta_2, ...$, $I'$ contains a regular sequence
consisting of forms $f_1,...,f_c$ of degrees $d'_1\le d_1\le
\delta_1, ..., d'_c\le d_c\le  \delta_c$.  Hence
$$B= \ell(R'/I') \le \ell(R'/(f_1,...,f_c)) =d'_1\cdots d'_c   \le \delta_1 \cdots  \delta_c.$$

b) Since $x_{c+1},...,x_{n}$ is a s.o.p. of $R/\Gin(I)$, it is
also a s.o.p. of $R/I$. Hence it is a minimal reduction of the
algebra $R/I$. By [V, Theorem 9.3.4]
$$x_i^{\adeg(I)} \in (I, x_{c+1},...,x_{n}) , \ \ \text{\rm for \ all} \ i \ge 1.$$
This means $x_1^{\adeg(I)},...,x_c^{\adeg(I)}$ form a regular
sequence in $I'$. The above argument gives $B \le (\adeg I)^c$.
\end{pf}

\begin{Remark}\label{A5} {\rm i) In the above theorem we may replace $b$ by
$(\reg I)^c$ in order to get  a bound for $h^j_S(t)$, which depends
only on $\reg I$ and $d,c$.  \par

ii)  Hilbert cohomological functions are of reverse polynomial
type, i.e. for each $i\ge 0$ there is a polynomial $p^i_S(t)$ such
that $h^i_S(t)= p^i_S(t)$ for all $t\ll 0$ (see [BrS, Theorem
17.1.9]). The number
$$\nu_S^i = \min\{ t\in {\Bbb Z};  \ h^i_S(t)\neq  p^i_S(t)\} -1$$
is called $i$-th cohomological postulation number of $S$ (see
[BrL]). Thus, if  $H^i_\mm(S)$, $i<d$ is of finite length, then all
graded components $H^i_\mm(S)_t$ vanish below $\nu_S^i $. Brodmann
and Lashgari proved that all $-\nu_S^i, \ i\le d, $ can be bounded
by a polynomial (of huge degree) in the numbers
$h^1_S(0),...,h^d_S(-d+1)$ (see [BrL, Theorem 4.6]). Combining
their result with Theorem \ref{A4} we see that $-\nu_S^i $ can be
bounded by a polynomial in $\reg S$. Thus, the number of
"irregular" negative components of local cohomology modules is
governed by the Castelnuovo-Mumford regularity. }\end{Remark}

\section{Hilbert coefficients} \label{B}\smallskip

Write the Hilbert polynomial in the  form:
$$P_{S}(t)= e_0{t+d-1 \choose d-1} -e_1 {t+d-2 \choose d-2}+ \cdots + (-1)^{d-1}e_{d-1}.$$
Then $e_0, e_1,...,e_{d-1}$ are called  {\it Hilbert coefficients}
of $S$. Note that $e_0=e$. Sometimes we also write $e_i = e_i(S)$
to emphasize its dependence on $S$.\par

We first  estimate $|e_i|$ in terms of the arithmetic degree. For
the application later, the following result is formulated in a
rather technical way.

\begin{Theorem}\label{B3} Let $K$ be an infinite field and $I$  an arbitrary homogeneous ideal. 
Assume that
$x_{c+1},...,x_{n}$ are chosen generically. Let $T_d = R/
(I:\mm^\infty),\ T_{d-1} = R/((I,x_n):\mm^\infty), ...$. Then
\begin{itemize}
\item[(i)] $|e_1| \le (\adeg I)^c (\reg T_2 +1) \leq (\adeg I)^c
\reg I$.
 \item[(ii)] For $i\ge 2$, $|e_i| \le \frac{3}{2}(\adeg
I)^c (\reg T_{i+1} +1)^i\le  \frac{3}{2}(\adeg I)^c ( \reg I)^i$.
\end{itemize}
\end{Theorem}

\begin{pf} The second inequalities in both (i) and (ii)  follow from Lemma \ref{B2}(i).
Let us prove the first ones. Set $T = T_d$. Let $H_T(t)$ denote
the Hilbert function of $T$. Since $S$ and $T$ have the same
Hilbert polynomial, $e_{i}= e_i(T)$ for all $i\ge 0$. From the
Grothendieck-Serre  formula
$$ P_T(t) - H_T(t) = \sum_{i=0}^d (-1)^{i+1}h^i_T(t),$$
we get (setting $t=-1$):
$$(-1)^{d-1}e_{d-1} = C - D,$$
where
$$C= h^1_T(-1) + h^3_T(-1) \ \cdots ,$$
and
$$D= h^2_T(-1) + h^4_T(-1) \ \cdots .$$
Hence
$$|e_{d-1}| \le \max\{C,\ D\}.$$
For short, set $b = (\adeg I)^c$. If $d=2$, then by Theorem
\ref{A4} we have $C = h^1_T(-1) \le (b-1)(\reg T +1)$ and $D=
h^2_T(-1) \le b\cdot \reg T$. Therefore $|e_1| \le b\cdot (\reg
T_2 +1)$.

Let $d\ge 3$. Since $x_n$ is generic, it is a regular element on
$T$. Since $P_{T_{d-1}}(t) = P_{T/x_nT}(t)$,  we have $e_i =
e_i(T) = e_i(T_{d-1})$ for all $i \le d-2$. The corresponding
sequence of rings constructed for $T_{d-1}$ as above are exactly
the rings $T_{d-1}, T_{d-2},...,T_1$. By Lemma \ref{B2}(ii),
$\adeg T_{d-1} \le \adeg T$. Hence, by the induction hypothesis, it
remains to prove (ii) for $i=d-1$. Note that
$${v+u-1 \choose u} \le v^u,\ \ \text{and} \ \ {v+1 \choose u} \le (v+1) \frac{v^{u-1}}{u!}.$$
Let $r = \reg T$. If $d= 2k+1$, where $k\ge 1$, then Theorem
\ref{A4} yields
$$\begin{array}{ll}
C &\le (b-1)(r+1)r^{d-2}\{ 1 + \frac{1}{3!} + \cdots +  \frac{1}{(2k-1)!}\} + b  
\frac{r^{d-1}}{(2k)!} \\
&\le b(r+1) r^{d-2}\{ 1 + \frac{1}{3!} + \cdots +  \frac{1}{(2k-1)!} + \frac{1}{(2k)!}\}\\
& \le \frac{3}{2} b(r+1)^{d-1},
 \end{array}$$
and
$$D \le (b-1)(r+1) r^{d-2}\{  \frac{1}{2!} + \cdots +  \frac{1}{(2k)!}\} \le 
(b-1)(r+1)^{d-1}.$$
Hence $|e_{d-1}| \le \frac{3}{2} b(r+1)^{d-1}$.

 The inequality in the case $d= 2k,\
k\ge 2$, can be shown similarly.
\end{pf}

\begin{Remark}\label{B5} {\rm a) In the above proof, if $h^1_T(-1)= 0$, then $C \le
(\adeg I)^c(r+1)^{d-1}$. Hence, if $h^1_{T_{i+1}} (-1) = 0$, then
$$|e_i| \le (\adeg I)^c (\reg T_{i+1} +1)^i\le  (\adeg I)^c ( \reg I)^i.$$

b) Consider again Example \ref{C9}: $S= K[x,y,u,v]/((x,y)^2,
xu^t+yv^t),\ t\ge 1$. We have $ e_1 = - (t+1)$, while $\reg (S) = t,\ \adeg S = 2$ and the bound in (i) 
of the above theorem is $4(t+1)$. Thus one cannot avoid $\reg I$ in the above theorem.
}\end{Remark}

Note that $\dim T_{i+1} = i+1$. Combining Theorem \ref{B3} and
Theorem \ref{C4} we get

\begin{Proposition}\label{C10} Let $S$ be a reduced ring of dimension at least two. Then
\begin{itemize}
\item[(i)] $|e_1| \le (\adeg S)^c(\frac{e(e-1)}{2} + \adeg S)$.
\item[(ii)] $|e_i| \le \frac{3}{2} (\adeg S)^c(\frac{e(e-1)}{2} +
\adeg S)^{i2^{i-1}}$ if $i\ge 2$.
\end{itemize}
\end{Proposition}

Another consequence of Theorem \ref{B3} is: 

\begin{Corollary}\label{B6} Assume that $K$ is an algebraically closed field of
characteristic 0 and $\text{\rm Proj}(R/I)$ is a reduced and
pure-dimensional smooth subscheme in ${\Bbb P}^{n-1}$. Then for
all $i\ge 1$ we have
$$ |e_i| < (i+2)^i e^{c+i}.$$
\end{Corollary}

\begin{pf} By Bertini's theorems (see [FOV, Corollary 3.4.6 and Corollary 3.4.14]) we may
 assume that all $\text{\rm Proj}(T_{i})$ are reduced and pure-dimensional
  smooth subschemes. By Mumford's bound: $\reg T_{i+1} \le (i+2)(e-2) + 1$
  (see [BM, Theorem 3.12(ii)]). Moreover, in this case $h^1_{ T_{i+1}}(-1) = 0$
  for all $i \ge 1$ and $\adeg I = e$ . Hence, by Theorem \ref{B3} and Remark \ref{B5}, we get
$$|e_i| \le e^c ((i+2)(e-2) + 2)^i < (i+2)^i e^{c+i}.$$
\end{pf}

\begin{Remark}\label{C11} {\rm Let $S$ be a reduced ring of dimension at least two.

i) It is known that for any $K$-algebra $S$, $e_1 \le e(e-1)/2$
(see [Bl, Remark 3.10]). Hence in the statement  (i) of
Proposition \ref{C10} only the following inequality is new: $e_1
\ge - (\adeg S)^c(\frac{e(e-1)}{2} + \adeg S)$.\par

ii) Let us recall\par

\noindent {\bf Eisenbud-Goto conjecture} [EG]: Let $K$ be an
algebraically closed field. If $I$ is a prime ideal containing no
linear form, then $\reg R/I \le e - c$.\par
 
If this conjecture holds true, then by Remark \ref{B5},   $|e_i| \le
 (\deg S)^{c+i}$ provided  $S$ is a domain. Note that the Eisenbud-Goto conjecture is close to be 
proved for smooth varieties of dimension at most 6 over a field of characteristic zero, by the work of 
several people including Lazarsfeld, Ran and Kwak.
 This   indicates that the bounds in
Theorem \ref{C4} and Proposition \ref{C10} are probably far from being
sharp.\par

iii) There is a bound on $|e_i|$ in terms of the so-called
homological degree which also holds for any standard graded
algebra over an artinian ring, see [RVV, Theorem 4.3]. However the
homological degree is very big.}\end{Remark}

We  now estimate $|e_i|$ by mean of the defining degrees. Recall
that homogeneous elements $y_1,...,y_m$ of $S$ form a {\it filter
regular sequence} if $[(y_1,...,y_{i-1}): y_i]_t =
(y_1,...,y_{i-1})_t$ for all $t \gg 0$ and $i=0,...,m$. On the
other words, $y_i$ is a filter regular element on
$S/(y_1,...,y_{i-1})S$.

\begin{Theorem}\label{B3b} Let $I$ be an arbitrary homogeneous ideal.
Assume that $d\ge 2$ and $x_{c+1},...,x_n$ is a filter regular
sequence on $S$. Let $S_d = S,\ S_{d-1} = S/x_nS_d,\ ... $  Set
$\pi =  \delta_1 \cdots  \delta_c$. Then
\begin{itemize}
\item[(i)] $|e_1| \le \pi\cdot (\reg S_2 +1) \leq \pi\cdot  \reg
I$. \item[(ii)] For $i\ge 2$, $|e_i| \le \frac{3}{2}\pi\cdot (\reg
S_{i+1} +1)^i\le  \frac{3}{2}\pi\cdot  ( \reg I)^i$.
\end{itemize}
\end{Theorem}

\begin{pf}
 The proof is similar to that of Theorem \ref{B3} after noticing that
$\reg S_i \le \reg S_{i+1}$ and $\delta_1(S_i) \le
\delta_1(S_{i+1}) ,..., \delta_c(S_i) \le \delta_c(S_{i+1})$ for
all $i \ge 2$.
\end{pf}

Combining it with Theorem \ref{E1}  we immediately get

\begin{Proposition}\label{E7}
Let $d\ge 1$. Then
\begin{itemize}
\item[(i)] $|e_1| \le  \pi(\pi + \Delta -1)^2$. 
\item[(ii)] For all
$i\ge 2$, we have $|e_i| \le \frac{3}{2}\pi (\pi + \Delta -1)^{i2^i}$.
\end{itemize}
In particular $|e_i| < (\frac{3}{2}\Delta^c + \Delta)^{1+i2^i}$ for all $i\ge 1$.
\end{Proposition}

 A direct application of Theorem \ref{B3b} sometimes  gives much better bounds than the ones in 
the previous proposition. For example, using [BEL]  and the second 
inequality in Theorem \ref{B3b}(ii), one immediately gets that 
$$|e_i| \le \frac{3}{2}\pi (\reg I)^i \le \frac{3}{2} \Delta^c(c(\Delta-1))^i <
\frac{3}{2} c^i\Delta^{c+i},$$
 provided  $\text{\rm Proj}(R/I)$ is a
reduced and pure-dimensional smooth subscheme. Another case is

\begin{Corollary}\label{B7} Let $I$ be an  ideal generated by monomials of degree at
most $\Delta$ in $n$ variables.  Then for all $i\ge 1$ we have
$$|e_i| \le \frac{3}{2} \min\{ (\adeg I)^{c+i}, \ n^i\Delta^{c+i} \}.$$
\end{Corollary}

\begin{pf} By [HT, Theorem 1.1], $\reg I \le \adeg I$ and by  Taylor's resolution (see also [HT, 
Theorem 1.2]), $\reg I \le n\Delta$.
 Hence the statement follows from theorems \ref{B3} and \ref{B3b}.
\end{pf}

The following example shows that the bounds in  Theorem \ref{B3b}
and Corollary \ref{B7} are rather good.

\begin{Example}\label{B4} {\rm Let $n> c+1$ and
$$I= (x_1,...,x_c) \cap (x_1^r,..., x_c^r, x_{c+1}^{r-1},...,x_{n-1}^{r-1}).$$
Using the exact sequence
$$0 \to R/I \to R/P \oplus R/J \to R/(P+J) \to 0,$$
 where $P = (x_1,...,x_c)$ and $J = (x_1^r,..., x_c^r, x_{c+1}^{r-1},...,x_{n-1}^{r-1})$,
 one can check that
$$\reg R/I = (n-1)r-2n+c+2,$$
and
$$P_{R/I}(t) = {t+d-1 \choose d-1} + [ r^c(r-1)^{d-1} - (r-1)^{d-1}].$$
Hence $|e_{d-1}| = (r^c-1)(r-1)^{d-1}$, while by Corollary
\ref{B7}  $|e_{d-1}| \le 3r^{c+d-1}n^{d-1}/2$, and by Theorem
\ref{B3b} $|e_{d-1}| \le 3 r^c[(n-1)r-2n+c+3]^{d-1}$/2.
}\end{Example}

 \section{Finiteness of Hilbert functions } \label{D}\smallskip

In this section we prove Theorem \ref{I1}. We need some further
preliminary results. The following result extends an estimation of $H_S(t)$ mentioned in Remark 
\ref{C6b} to arbitrary ideals.

\begin{Lemma}\label{D1} Let $I$ be an arbitrary homogeneous ideal. Let
$$b = \min \{ \delta_1\cdots \delta_c,\ (\adeg S)^c\}.$$
For all $t\ge 0$ we have
$$H_S(t) \le (b-1){t+d-2\choose d-1} + {t+d-1\choose d-1}.$$
\end{Lemma}

\begin{pf} We may assume that $x_{c+1},...,x_n$ are chosen generically. In particular,
  $x_{c+1},...,x_n$ form a s.o.p. of $S$. Set
 $B = \ell(S/(x_{c+1},...,x_n)S)$. By [RVV, Proposition 2.4] for all $t\ge 0$ we have
$$H_S(t) \le (B -1){t+d-2\choose d-1} + {t+d-1\choose d-1}.$$
As shown in the proof of Theorem \ref{A4}, $B\le b$. Hence the
lemma is proven.
\end{pf}

\begin{Lemma}\label{D2} Assume that $K$ is an algebraically closed field, $I$
is an intersection of prime ideals and $I$ contains no linear
form. Then $c \le d(\adeg I -1)$.
\end{Lemma}

\begin{pf} By the assumption
$$ I = \cap_{i=1}^s \pp_i,$$
where $\pp_i$ are prime ideals of height at least $c$. Since $s
\le \adeg I$, the statement is derived from the following
inequality
$$c \le \adeg I - s + (s-1)d.$$
We prove this inequality by induction on $s$. The case $s=1$ is
well known. Let $s>1$. Put $J = \cap_{i=1}^{s-1} \pp_i.$ Let $a$
and $b$ be the maximal number of independent linear forms
contained in $J$ and $\pp_s$, respectively. By the induction hypothesis, we have
$c- a \le \adeg J - (s-1) + (s-2)d$ and $ c-b \le e(R/\pp_s)-1$.
Since $\adeg I = \adeg J + e(R/\pp_s)$, we get
$$2c \le \adeg I -s + (s-2)d + a+b.$$
If $a+b> n$ it would imply that there is a linear form in $J \cap
\pp_s = I$, a contradiction. Hence $a+b\le n = d+c$. The above
inequality then yields $c \le \adeg I - s + (s-1)d.$
\end{pf}

As  mentioned in the introduction, Theorem \ref{I1} is an
immediate consequence of Theorem \ref{C4}, Lemma \ref{D2} and
[RTV2, Theorem 2.3]. We give here a direct proof without the use
of [RTV2].\vskip0.5cm

\noindent {\it Proof of Theorem \ref{I1}.} Without the loss of
generality we may assume from the beginning that $I$ contains no
linear form. Note that $e \le \adeg S$. Therefore, by Proposition
\ref{C10} and Lemma \ref{D2}, there are only finitely many Hilbert
polynomials associated to reduced algebras such that $\adeg S \le
a$ and $\dim S \le d$. By Lemmas \ref{D1}  and \ref{D2}, there are
only finitely many choices for the initial values of  Hilbert
functions, while Theorem \ref{C4} says that for $t \ge
(\frac{e(e-1)}{2} + \adeg S)^{2^{d-2}}$ each Hilbert function
agrees with the corresponding Hilbert polynomial. This implies the
finiteness of the number of Hilbert functions. \hfill $\square$.
\vskip0.5cm

Example \ref{C9} shows that without the assumption $S$ being a
reduced ring Theorem \ref{I1} does not hold.

Applying Proposition \ref{E7} and Theorem \ref{E1}, as in the
proof of Theorem \ref{I1}, we  get a similar finiteness result in
terms of the defining degrees.

\begin{Corollary}\label{E8} Given two numbers $\delta$ and $n$, there exist
only finitely many Hilbert functions associated to homogeneous ideals generated
 by forms of degrees at most $\delta$ in at most $n$ variables.
\end{Corollary}

\section{Castelnuovo-Mumford regularity of initial ideals} \label{F}\smallskip

In this last section we apply results in the previous sections to
study the Castelnuovo-Mumford regularity of an initial ideal $\gin
(I)$ of $I$ with respect to any given term order and coordinates. We
even consider a much bigger class: the class of all ideals $J$
having the same Hilbert function as  $I$. Then one can easily
bound $\reg J$ in terms of some data of $I$. This approach was
initiated in [CM] and developed further in [HH]. Let us recall
some notations. The Hilbert polynomial can be uniquely written in
the form
$$P_{R/I}(t)= {c_1+ t \choose t} + {c_2+ t-1 \choose t-1}+ \cdots + {c_s+ t-s+1 \choose t-s+1},$$
where $c_1\geq c_2 \ge \cdots \ge c_s \geq 0$ are integers (see,
e.g., [V, Section B6]).  For $0\le i \le d-1$ set
$$B_i = \sharp\{j;\  c_j \ge (d-1)- i\}.$$
Thus in the above notations, $s= B_{d-1}$ (for convenience, we set
$B_{-1} = 0$). The following result easily follows from Gotzmann's
regularity theorem:

\begin{Lemma}\label{F1} {\rm [HH, Lemma 5]} Let $I,J$ be homogeneous ideals having
 the same Hilbert function. Then
$$\reg J \le  \max\{ \reg I,\ B_{d-1}\}.$$
\end{Lemma}

Since we already know bounds for $\reg(I)$ (see Theorems \ref{C4}
and \ref{E1}), we have only to estimate $B_{d-1}$.  For this
purpose we need some  relations between the invariants $B_i$ just
defined and the Hilbert coefficients which were given in  [Bl,
Proposition 3.9] (see also [CM, Lemma 1.5]).

\begin{Lemma}\label{F2}   For all $0 \le j\le d-1$ we have
$$B_i = (-1)^ie_i + {B_{i-1}+1 \choose 2} - {B_{i-2}+1 \choose 3} + \cdots +
(-1)^{i+1}{B_{0}+1 \choose i+1}.$$
\end{Lemma}

Note that  $B_{d-1} \ge \cdots \ge B_0 = e$. In order to estimate
$B_j$, we need the following combinatorial result.

\begin{Lemma}\label{F3} Assume that
$$|e_i| \le M^{\alpha + i\beta 2^i} \ \ \text{for\ all\ } i\ge 0,$$
where $M\ge 2$ and $\alpha,  \beta \ge 1$. Then for all $0 \le j
\le d-1$ we have
$$B_j \le M^{(\alpha+j\beta)2^j}.$$
\end{Lemma}

\begin{pf} We have $B_0 = e = e_0 \le M^\alpha$ by the assumption. By Lemma 
\ref{F2}
the following holds
$$B_1 = -e_1 + {B_0 +1 \choose 2} \le  |e_1| + \frac{e(e+1)}{2}
< M^{\alpha + 2\beta} +  M^{2\alpha} \le M^{2(\alpha + \beta)}\ \
(\text{since}\ M\ge 2).$$

Let $j\ge 2$. Assume that
\begin{eqnarray} \label{EF1} B_{j-l} \le M^{(\alpha+(j-l)\beta)2^{j-l}}
\end{eqnarray}
for all $l\ge 1$. Lemma \ref{F2} yields:
\begin{eqnarray}
B_j &=& (-1)^je_j + {B_{j-1}+1 \choose 2} - {B_{j-2}+1 \choose 3} + \cdots +
(-1)^{j-1}{B_{0}+1 \choose j-1} \nonumber\\
& \le & |e_j| + {B_{j-1}+1 \choose 2} + {B_{j-3}+1 \choose 4} +
\cdots \label{EF2}
\end{eqnarray}
By (\ref{EF1}), for all $1\le l\le j$ we have
\begin{eqnarray} \label{EF4}
{B_{j-l}+1 \choose l+1} \le  \frac{(B_{j-l}+1)^{l+1}}{(l+1)!} \le
\frac{(B_{j-l}+1)^{2^l}}{(l+1)!} \le
\frac{M^{(\alpha+j\beta)2^j}}{(l+1)!}.
\end{eqnarray}
From (\ref{EF2}),  (\ref{EF4})  and the assumption $|e_j| \le
M^{\alpha + j\beta 2^j}$ it follows that
$$\begin{array}{ll}
B_j & \le M^{\alpha +j\beta2^j} + M^{(\alpha+j\beta)2^j} \{ \frac{1}{2!} +
 \frac{1}{4!} \ \cdots \}\\
&< M^{\alpha +j\beta2^j} + \frac{2}{3} M^{(\alpha+j\beta)2^j}\\
& \le M^{(\alpha+j\beta)2^j}.
 \end{array}$$
\end{pf}

By Macaulay's theorem $H_{R/\gin I}(t) =  H_{R/ I}(t)$ for all $t
\in \Bbb Z$. Hence, Theorem \ref{I2} stated in the introduction is
a special case of the following result.

\begin{Theorem}\label{F4} Let $K$ be an arbitrary field. Let $J$ be an arbitrary
homogeneous ideal of $ R=K[x_1,...,x_n]$ such that $H_{R/J}(t) = H_{R/I}(t)$ for all  $t$. Then
\begin{itemize}
\item[(i)]  $\reg(J) \le  (\frac{3}{2} \Delta^c + \Delta)^{d2^{d-1}}.$ \item[(ii)]
Moreover, if $R/I$ is a reduced algebra, then we also have
$$\reg(J) \le (\adeg(I))^{(n-1)2^{d-1}}.$$
\end{itemize}
\end{Theorem}

\begin{pf} 
(i) By Proposition \ref{E7}, $|e_i| < (\frac{3}{2}\Delta^c + \Delta)^{1+i2^i}$ for
all $i\ge 0$. Applying Lemma \ref{F3} to $M= \frac{3}{2}\Delta^c + \Delta,\ \alpha
=1,\  \beta = 1$ and $j=d-1$, we get $B_{d-1}\le
(\frac{3}{2}\Delta^c + \Delta)^{d2^{d-1}}.$ Then (i) follows from Lemma \ref{F1} and
Theorem \ref{E1}.

(ii) For short, set $a = \adeg I$. Note that  $a \ge e$ and
\begin{eqnarray} \label{EF5}
 \frac{e(e-1)}{2} + a \le   a^2.
\end{eqnarray}
Hence, by Proposition \ref{C10}(i)
$$|e_1| \le a ^{c+2}.$$
Let $i \ge 2$. Since $\Delta\ge 2$,  $a \ge 2$. By Proposition
\ref{C10}(ii) and (\ref{EF5}), we have
$$\begin{array}{ll}
|e_i| &\le \frac{3}{2}a^c (\frac{a(a+1)}{2})^{i2^{i-1}} = a^c
(\frac{a(a+1)}{2})^{i2^{i-1}-4}\cdot [\frac{3}{2}(\frac{a(a+1)}{2})^4]\\
& \le a^{c + 2(i2^{i-1}-4)} a^8 = a^{c+i2^i}.
 \end{array}$$
Thus, applying Lemma \ref{F3} to $M=a,\ \alpha = c,\ \beta =1$ and
$j = d-1$, we get $B_{d-1} \le a^{(n-1)2^{d-1}}$. By Lemma \ref{F1}
and Theorem \ref{C4} this implies $\reg J\le a^{(n-1)2^{d-1}}$.
\end{pf}

Note that if $R/I$ is a Cohen-Macaulay ring of dimension $d\ge 2$
(but not necessarily reduced), then one can get a little bit
better bound (see [HH, Theorem 9]):
$$\reg J \le e^{2^{d-1}}/ 2^{2^{d-2}}.$$

\begin{Example}\label{F5} {\rm Let $I^{lex} $ denote the lex-segment ideal
associated to the Hilbert function $H_{R/I}(t)$. This is the ideal
generated by all first $H_{I}(m)$ monomials of degrees $m$ with
respect to the lexicographic order, when $m$ runs through all
positive integers. It has the same Hilbert function as $I$. If
$R/I$ is a Cohen-Macaulay ring of dimension $d\ge 2$, then from
[CM, Theorem 2.5] it follows that $\reg (I^{lex} ) = B_{d-1}$.\par

i) Let $I$ be an ideal  generated by a regular sequence consisting
of forms of degrees $\delta_1 \ge \cdots \ge \delta_c$ such that
$c\ge 2$ and $\delta_2\ge 35$ ($d\ge 2$). It was shown in [HH,
Example 13] that
$$\reg (I^{lex} )  \ge 9\frac{\Delta^{c2^{d-1}}}{9^{2^{d-2}}}.$$
This shows that the bound in Theorem \ref{F4}(i) is close to be sharp.

ii) Let $S$ be a Veronesian embedding $K[y_1,...,y_d]^{(p)}$, i.e.
$S_1$ is generated by all monomials of degree $p$ in the variables
$y_1,...,y_d$, where $ d\ge 3$. This is a Cohen-Macaulay domain
and $P_S(t) = {pt +d-1 \choose d-1}$. Hence $\adeg S = e= p^{d-1}$
and $e_1 = dp^{d-2}(p-1)$. Let $p\ge 35$. Then $e_1 < e^2/36$ and
$e\ge 35^2$. Let $S= K[x_1,...,x_q]/I$, where $q = {p+d-1\choose
d-1}$. By [HH, Proposition 12] we get
$$\reg (I^{lex} )  \ge 9\frac{e^{2^{d-1}}}{9^{2^{d-2}}}.$$
This shows that the bound in the second part of Theorem \ref{F4}
is  close to be sharp too. }
\end{Example}

Since $\reg (\gin I) \ge \reg I$, the ideals of Mayr and Meyer
again show that the bound $(2\Delta^c)^{d2^{d-1}}$ of Theorem
\ref{I2} is rather good (see Remark \ref{E6}). We do not know
whether one can construct a reduced algebra $R/I$ such that there
is a term order with $\reg(\gin I)$ close to
$(\adeg(I))^{(n-1)2^{d-1}}$.

Finally we would like to make the following remark: In the proof
of Theorem \ref{I2} we use very rough estimation for $\reg I$ and
$|e_i|$. It could suspect that if $\reg I$ and $|e_i|$ are small,
then one could get a bound for $\reg(\gin I)$, which would be a
single exponent of $d$. But this is not the case as shown by [HH,
Section 4].

 \section*{References}\smallskip

\begin{itemize}

\item[[BM]] D. Bayer and D. Mumford, {\it  What can be computed in
algebraic geometry?}  Computational algebraic geometry and
commutative algebra (Cortona, 1991),  1--48, Sympos. Math., XXXIV,
Cambridge Univ. Press, Cambridge, 1993; MR 95d:13032.

\item[[BS]]  D. Bayer and M. Stillman, {\it A criterion for
detecting $m$-regularity}, Invent. Math. {\bf 87}(1987), no. 1,
1--11, MR 87k:13019.

\item[[BEL]]  A. Bertram, L. Ein and R. Lazarsfeld,  {\it Vanishing
theorems, a theorem of Severi, and the equations defining
projective varieties}, J. Amer. Math. Soc. {\bf 4}(1991),
587--602; MR 92g:14014.

 \item[[Bl]] C. Blancafort,  {\it Hilbert functions of
graded algebras over Artinian rings},  J. Pure Appl. Algebra  {\bf
125}(1998),  no. 1-3, 55--78; MR 98m:13023.

\item[[BrL]]   M. P. Brodmann and A. F. Lashgari,  {\it A diagonal
bound for cohomological postulation numbers of projective
schemes,}  J. Algebra  {\bf 265}(2003), 631--650; MR 2004f:14030.

\item[[BrS]]   M. P. Brodmann and R. Y. Sharp,   Local cohomology:
an algebraic introduction with geometric applications. Cambridge
Studies in Advanced Mathematics, 60. Cambridge University Press,
Cambridge, 1998; MR 99h:13020.

\item[[CS]] G. Caviglia and E. Sbarra, {\it  Characteristic-free
bounds for Castelnuovo-Mumford regularity},   Compos. Math. 141 (2005), no. 6, 1365--1373; MR 
2006i:13032.

\item[[CF]] M. Chardin and A. L. Fall, {\it Sur la r\'egularit\'e de Castelnuovo-Mumford de id\'eaux 
en dimension 2},    C. R. Math. Acad. Sci. Paris 341 (2005), no. 4, 233--238; MR 2164678.

\item[[CM]] M. Chardin and G. Moreno-Socias,  {\it  Regularity of
lex-segment ideals: some closed formulas and applications},
Proc. Amer. Math. Soc. {\bf 131}(2003),  no. 4, 1093--1102; MR 2003m:13014.

\item[[E]]  D. Eisenbud,  Commutative algebra. With a view toward
algebraic geometry. Graduate Texts in Mathematics, 150.
Springer-Verlag, New York, 1995; MR 97a:13001.

\item[[EG]] D. Eisenbud and S. Goto,  {\it Linear free resolutions
and minimal multiplicity}, J. Algebra {\bf 88}(1984), no. 1,
89--133; MR 85f:13023.

\item[[FOV]] H. Flenner, L. O'Carroll and W. Vogel,  Joins and
intersections. Springer Monographs in Mathematics.
Springer-Verlag, Berlin, 1999; MR 2001b:14010.

\item[[Gi]]  M. Giusti,  {\it Some effectivity problems in
polynomial ideal theory.}  EUROSAM 84 (Cambridge, 1984), 159--171,
Lecture Notes in Comput. Sci.,  {\bf 174}, Springer, Berlin, 1984,
MR 86d:12001.

\item[[GLP]] L. Gruson, R. Lazarsfeld and C. Peskine,  {\it On a
theorem of Castelnuovo, and the equations defining space curves},
Invent. Math.  {\bf 72}(1983),  no. 3, 491--506; MR 85g:14033.

\item[[HPV]]    J. Herzog, D.  Popescu and M.  Vladoiu,   {\it   On
the Ext-modules of ideals of Borel type}, Commutative algebra
(Grenoble/Lyon, 2001), 171--186, Contemp. Math., 331, Amer. Math.
Soc., Providence, RI, 2003; MR 2013165.

\item[[HH]]  L. T.  Hoa and E. Hyry, {\it Castelnuovo-Mumford
regularity of initial ideals}, J. Symb. Comp. {\bf 38}(2004), 1327--1341; MR 2168718.

\item[[HSV]]   L. T. Hoa, J. St\"uckrad  and W.  Vogel, {\it Towards
a structure theory for projective varieties of
 degree =codimension +2}, J. Pure Appl. Algebra,
{\bf 71}, 203-231(1991). MR 92f:14002.

\item[[HT]]  L. T. Hoa; N. V. Trung, {\it  On the
Castelnuovo-Mumford regularity and the arithmetic degree of
monomial ideals,}  Math. Z.  {\bf 229}(1998),  no. 3, 519--537; MR
99k:13034.

\item[[K]]  S. L. Kleiman,  {\it Les th\'eor\`emes de finitude
pour le foncteur de Picard} , in: "Th\'eorie des intersections et
th\'eor\`eme de Riemann-Roch" (French) S\'eminaire de
G\'eom\'etrie Alg\'ebrique du Bois-Marie 1966--1967 (SGA 6), pp.
616--666.  Lecture Notes in Mathematics, Vol. {\bf 225}.
Springer-Verlag, Berlin-New York, 1971; MR 50 $\sharp$7133.

\item[[MVY]]  C. Miyazaki, W. Vogel and K.  Yanagawa,  {\it
Associated primes and arithmetic degrees},  J. Algebra  {\bf
192}(1997),  no. 1, 166--182; MR 98i:13036.

\item[[MM]] H. M. M\"oller and F. Mora, {\it Upper and lower bounds for the degree of  Gr\"obner 
bases}, EUROSAM 84 (Cambridge, 1984), 172--183,
Lecture Notes in Comput. Sci.,  {\bf 174}, Springer, Berlin, 1984;
MR 86k:13008.

\item[[M]] D. Mumford, Lectures on curves on an algebraic surfaces, Princeton Univ. Press, 
Princeton 1966; MR 35$\sharp$187.

\item[[RTV1]]  M. E. Rossi, N. V. Trung and G. Valla, {\it
Castelnuovo-Mumford regularity and extended degree,}  Trans. Math.
Amer. Soc.  {\bf 355}(2003),  no. 5, 1773--1786; MR 2004b:13020.

\item[[RTV2]]  M. E. Rossi, N. V. Trung and G. Valla, {\it
Castelnuovo-Mumford regularity and finiteness of Hilbert
functions,}  Lect. Notes Pure Appl. Math., 244, Chapman \& Hall/CRC, Boca Raton, FL, 2006;  MR 
2184798.

\item[[RVV]]  M. E. Rossi, G. Valla and W. V. Vasconcelos, {\it
Maximal Hilbert functions,}  Results Math.  {\bf 39}(2001),  no.
1-2, 99--114; MR 2001m:13020.

\item[[S]]   E.  Sbarra,   {\it   Upper bounds for local
cohomology for rings with given Hilbert function},   Comm. Algebra
{\bf  29}(2001), no. 12, 5383--5409; MR 2002j:13024.

\item[[Sj]]  R. Sj\"ogren, {\it On the regularity of graded $k$-algebras of Krull dimension $\le 1$}, 
Math. Scand. {\bf 71}(1992), 167--172; MR 94b:13010.

\item[[V]]  W. V. Vasconcelos,  Computational methods in
commutative algebra and algebraic geometry. With chapters by D.
Eisenbud, D. R. Grayson, J. Herzog and M. Stillman. Algorithms and
Computation in Mathematics, 2. Springer-Verlag, Berlin, 1998; MR
99c:13048.

\end{itemize}
\end{document}